\documentclass[preprint,showpacs,showkeywords,preprintnumbers,amsmath,amssymb]{revtex4}

\usepackage{graphicx}
\usepackage{dcolumn}
\usepackage{bm}
\begin{document}

\title{Critical Properties of $S^{4}$ System Restudied via Generalized Migdal-Kadanoff Bond-moving Renormalization}

\author{Chun-Yang Wang£¨Íõ´ºÑô£©\footnote{Corresponding author. Electronic mail:
wchy@mail.bnu.edu.cn}}
\author {Wen-Xian Yang£¨ÑîÎÄÏ×£©}
\author {Hong Du£¨¶Åºì£©}

\affiliation{Shandong Provincial Key Laboratory of Laser
Polarization and Information Technology, College of Physics and
Engineering, Qufu Normal University, Qufu 273165, China}


\begin{abstract}
We study the critical properties of the spin-continuous $S^{4}$
system on the typical translational invariant triangular lattices by
combining the recently-developed generalized Migdal-Kadanoff
bond-moving recursion procedures with the cumulative expansion
technique. In three different cases of nearest-neighbor, next
nearest neighbor and external field we obtain the critical points
and further calculate the critical exponents according to the
scaling theory. In all case it is found that there exists three
fixed points. The correlation length critical exponents obtained
near the Wilson-Fisher fixed points are found getting smaller and
smaller with the increasing of the system complexity. Others are
found similar to the results of the classical Gaussian model and
different from those of the Ising system.
\end{abstract}



\maketitle

\section{INTRODUCTION}
The $S^{4}$ model, which is another extension of the classical Ising
model, is a kind of typical spin-continuous systems which allows the
spin to take any real value between $(-\infty,+\infty)$ instead of a
discrete one. The study of it is of great value for better
understanding the properties of the ferromagnetic systems from
either a theoretical or a practical viewpoint. Due to this point its
it has been widely studied in the past few years
\cite{msrg,liy,yin}. For an example, critical properties of the
$S^{4}$ model on some lattice systems with translational symmetry
has been found to have close dependent on the space dimensionality
$d$ : when $d$ is more than 4, only the Gaussian fixed point is
obtained while for low-dimensional systems with $d$ less equal to 4,
not only the Wilson-Fisher fixed point but also the Gaussian one can
be obtained \cite{msrg}. Similar dimensionality dependent behaviors
are also found for the $S^{4}$ model on the fractal lattices
\cite{liy}.

However, critical behaviors of the $S^{4}$ system is far from having
been well studied. This may in some case due to the self-complexity
of the model caused by its four-body interactions. Other
restrictions may come from the particular lattice system on which
the model is constructed. Recently, a type of generalization of the
remarkable Migdal-Kadanoff bond-moving renormalization procedure is
developed which enables us to handle complex systems in a relatively
easy way \cite{chyw3}. For the great convenience it has brought in
the coarse-graining process, in this paper, we extend this
renormalization procedure to study the critical properties of the
$S^{4}$ system on translational invariant triangular lattice by
combining it with the cumulative expansion method
\cite{liy,yin,liu}.

The paper is organized as follows: in Sec. \ref{sec2}, a brief
review of the generalized bond-moving recursion procedures is
presented; in Sec. \ref{sec3}, the procedures are used to study the
critical properties of the spin-continuous $S^{4}$ model constructed
on the translational invariant triangular lattice; Sec. \ref{sec4}
serves as a summary of our conclusion in which further implicit
applications are also discussed.

\section{renormalization procedure and cumulative expansion}\label{sec2}

\begin{figure}
\includegraphics[scale=0.8]{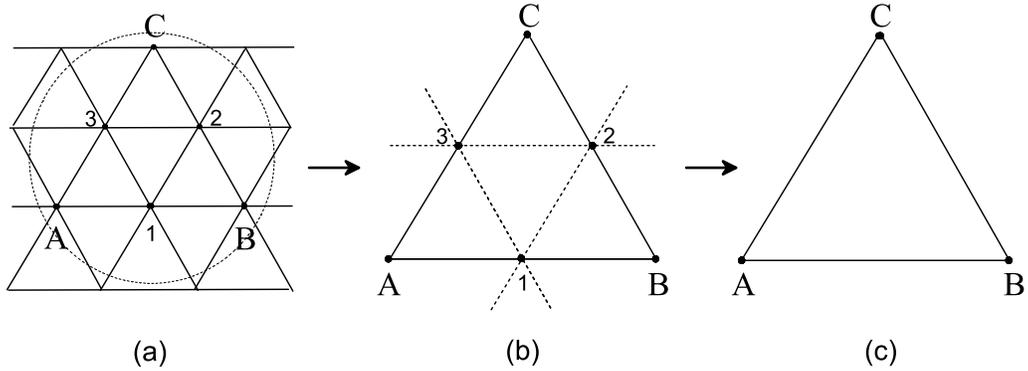}
\caption{Generalized bond-moving procedure recurring on the
triangular lattice where the peripheral bonds connecting the three
to be eliminated sites 1, 2 and 3 around the selected triplet
$\triangle ABC$ are drawn at half length. \label{pro1}}
\end{figure}

In Ref.\cite{chyw3}, a new type of generalization of the remarkable
Migdal-Kadanoff bond-moving renormalization procedure is developed
bringing with much convenience for the study of spin-continues
systems. By introducing in some kind of symmetrical half-length bond
operations, the renormalization process is greatly simplified as has
been witnessed in the study of Gaussian system. For the example of
systems constructed on the triangular lattice as illustrated in
Fig.\ref{pro1}, this was realized by supposing each bond connecting
the to be eliminated sites denoted as 1, 2 and 3 in $\triangle ABC$
to schlep symmetrically two half-length bonds. Thus the number of
lattice sites that a basic unit for recursion contains can be
reduced to only six which provides a good foundation for the
following integral calculation.

The technique of cumulative expansion is always used as an auxiliary
method in calculating high order integrations \cite{liy,yin,liu}.
For a system with effective Hamiltonian: $H_{\textrm{eff}}=H_{0}+V$,
a partial trace (PT) is defined as
\begin{eqnarray}
(PT)&=&\int^{+\infty}_{-\infty}\left(\prod^{m}_{i}ds_{i}\right)e^{H_{\textrm{eff}}}\\\nonumber
&=&\int^{+\infty}_{-\infty}\left(\prod^{m}_{i}ds_{i}\right)e^{H_{0}}
\frac{\int^{+\infty}_{-\infty}\left(\prod^{m}_{i}ds_{i}\right)e^{H_{0}+V}}
{\int^{+\infty}_{-\infty}\left(\prod^{m}_{i}ds_{i}\right)e^{H_{0}}}\\\nonumber
&=&A\langle e^{V}\rangle_{0},\label{pt}
\end{eqnarray}
where
\begin{eqnarray}
\langle\cdots\rangle_{0}=\frac{\int^{+\infty}_{-\infty}\left(\prod^{m}_{i}ds_{i}\right)(\cdots)e^{H_{0}}}
{\int^{+\infty}_{-\infty}\left(\prod^{m}_{i}ds_{i}\right)e^{H_{0}}}\label{pt2}
\end{eqnarray}
is namely the cumulative expansion average and
\begin{eqnarray}
A=\int^{+\infty}_{-\infty}\left(\prod^{m}_{i}ds_{i}\right)e^{H_{0}}.\label{pt0}
\end{eqnarray}
The integrations in the above equations trace over all the to be
decimated sites in the renormalization process.

On the supposition of $V$ is a small variable, we can expand $e^{V}$
by a Maclaurin series as
\begin{eqnarray}
e^{V}=1+V+\frac{1}{2!}V^{2}+\frac{1}{3!}V^{3}+\cdots.\label{pte}
\end{eqnarray}
The partial trace is then obtained to be
\begin{eqnarray}
(PT)=A\left(1+\langle V\rangle_{0}+\frac{1}{2!}\langle
V^{2}\rangle_{0}+\frac{1}{3!}\langle
V^{3}\rangle_{0}+\cdots\right).\label{pt3}
\end{eqnarray}
Noticing that the partition function should be kept unchanged after
the renormalization group (RG) transformation. That is, supposing
\begin{eqnarray}
\int^{+\infty}_{-\infty}\left(\prod^{m}_{i}ds_{i}\right)e^{H_{\textrm{eff}}}=Ce^{H'_{\textrm{eff}}},\label{pte}
\end{eqnarray}
then $H'_{\textrm{eff}}$ represents an effective Hamiltonian after
RG transformation ($C$ is a RG constant). Thus we can obtain from
the above equations
\begin{eqnarray}
H'_{\textrm{eff}}=\textrm{ln}A+\langle
V\rangle_{0}+\frac{1}{2}\left(\langle V^{2}\rangle_{0}-\langle
V\rangle^{2}_{0}\right)+\cdots,\label{ptee}
\end{eqnarray}
where the approximation relation
$\textrm{ln}(1+x)=x-\frac{x^{2}}{2}+\frac{x^{3}}{3}-\cdots$ is used.

Now we have found a convenient way to calculate the effective
Hamiltonian $H'_{\textrm{eff}}$ from the partial trace defined
herein before. Supposing
\begin{eqnarray}
\int^{+\infty}_{-\infty}\left(\prod^{m}_{i}ds_{i}\right)e^{H_{0}}=e^{H'_{0}},\label{pt00}
\end{eqnarray}
the first part in the right-hand site of Eq.(\ref{ptee}) can be
regarded as the zero order approximation of $H'_{\textrm{eff}}$.
Similar understandings are still true as well for other right-hand
parts. In the following sections, this will be applied to study the
critical behavior of the $S^{4}$ system constructed on the
triangular lattice combining with the generalized bond-moving
renormalization procedures.

\section{Critical behavior of the $S^{4}$ system} \label{sec3}

The $S^{4}$ model is another spin-continuous extension of the Ising
model \cite{s4}. By introducing a probability distribution function
\begin{eqnarray}
W(s_{1},s_{2},\cdots,s_{N})=e^{-\frac{b}{2}\sum_{i}s_{i}^{2}-\sum_{i}us_{i}^{4}},\label{pd}
\end{eqnarray}
for the classical Ising spin system, the effective Hamiltonian of
the $S^{4}$ model in an external field $h$ can be obtained as
\begin{eqnarray}
H_{\textrm{eff}}=\sum_{\langle
ij\rangle}Ks_{i}s_{j}-\frac{b}{2}\sum_{i}s_{i}^{2}-u\sum_{i}s_{i}^{4}+h\sum_{i}s_{i},\label{eH}
\end{eqnarray}
where $K={J}/{k_{B}T}$ is the reduced interaction with $K>0$ denotes
the ferromagnetic system; $k_{B}$ is the Boltzmann constant and $T$
the thermodynamic temperature; $h$ is the external field. The
$\langle ij\rangle$ in the first summation generally represents a
certain nearest-neighbor spin pair. Seen from Eq.(\ref{eH}), the
$S^{4}$ model Hamiltonian contains synchronously the four-spins and
two-spins interactions denoted by $u$ $(u>0)$ and Gaussian constant
$b$ respectively. This will undoubtedly made it more difficult to be
studied, however these can all be settled following the generalized
bond-moving renormalization procedures.

\subsection{nearest-neighbor interaction} \label{sbsec31}

\begin{figure}
\centering
\includegraphics[scale=1.2]{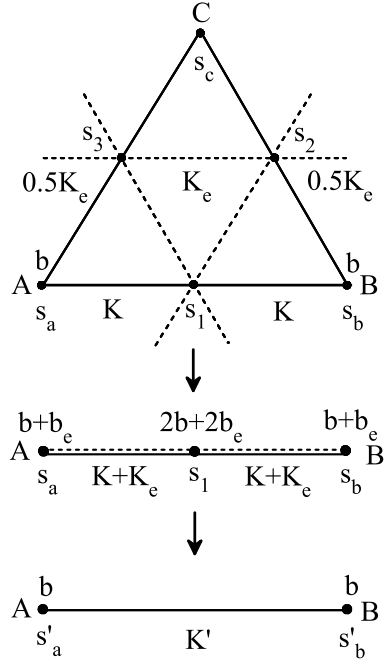}
\caption{Bond-moving and decimation procedures for the renormalized
bond $K'$ between sites A and B on the triangular
lattices.\label{pro2}}
\end{figure}

Firstly we give our study on the case of only considering
nearest-neighbor interactions. The procedure of the generalized
bond-moving RG transformation for obtaining a certain bond $K'$ is
shown in Fig.\ref{pro2} as an example.
From which we can see that the number of to be decimated sites
reduces to be only one. Thus the effective Hamiltonian of the basic
unit for RG transformation in the case of $h=0$ can be given as
\begin{eqnarray}
H_{\textrm{eff}}=2K\left(s_{a}s_{1}+s_{1}s_{b}\right)
-b\left(2s^{2}_{1}+s^{2}_{a}+s^{2}_{b}\right)
-2u\left(2s^{4}_{1}+s^{4}_{a}+s^{4}_{b}\right),\label{eeH}
\end{eqnarray}
where in the RG procedures two types of interactions $K_{e}$ and $K$
together with two types of self-energy $(-b_{e}s^{2}/2)$ and
$(-bs^{2}/2)$ are assigned respectively for differentiation of the
to be and not to be eliminated bonds Gefen et al did on the fractal
\cite{Genfen,Genfen2,Genfen3}. For the particular case of triangular
lattices the numerical value of $K_{e}$ and $K$ is actually
identical as well as that of $b_{e}$ and $b$ (or $u_{e}$ and $u$).
In the recursion procedures the two half length bonds are considered
acting effectively as a whole one and be moved regularly as other
bonds.

Supposing
\begin{eqnarray}
V=-4us^{4}_{1},\label{V}
\end{eqnarray}
and
\begin{eqnarray}
H_{0}=2K\left(s_{a}s_{1}+s_{1}s_{b}\right)
-b\left(2s^{2}_{1}+s^{2}_{a}+s^{2}_{b}\right)
-2u\left(s^{4}_{a}+s^{4}_{b}\right),\label{H0}
\end{eqnarray}
the zero-order approximation of the effective Hamiltonian after RG
transformation can be obtained from Eq.(\ref{pt00}) as
\begin{eqnarray}
H'_{0}=K_{01}s_{a}s_{b}+K_{02}\left(s^{2}_{a}+s^{2}_{b}\right)
+K_{03}\left(s^{4}_{a}+s^{4}_{b}\right),\label{H'0}
\end{eqnarray}
with
\begin{subequations}\begin{eqnarray}
K_{01}&=&\frac{K^{2}}{b},\\
K_{02}&=&\frac{K^{2}}{2b}-b,\\
K_{03}&=&-2u.
\end{eqnarray}\label{00}\end{subequations}
Continuously we can obtain the first-order one by Eq.(\ref{pt2}) as
\begin{eqnarray}
\langle
V\rangle_{0}&=&\frac{\int^{+\infty}_{-\infty}\left(\prod^{m}_{i}ds_{i}\right)Ve^{H_{0}}}
{\int^{+\infty}_{-\infty}\left(\prod^{m}_{i}ds_{i}\right)e^{H_{0}}}\\\nonumber
&=&K_{11}s_{a}s_{b}+K_{12}\left(s^{2}_{a}+s^{2}_{b}\right)
+K_{13}\left(s^{4}_{a}+s^{4}_{b}\right),\label{H'1}
\end{eqnarray}
where
\begin{subequations}\begin{eqnarray}
K_{11}&=&-\frac{3K^{2}u}{b^{3}},\\
K_{12}&=&-\frac{3K^{2}u}{2b^{3}},\\
K_{13}&=&-\frac{K^{4}u}{4b^{4}}.
\end{eqnarray}\label{10}\end{subequations}
In the same way, the second-order term can also be approximately
obtained as
\begin{eqnarray}
\frac{1}{2}\left(\langle V^{2}\rangle_{0}-\langle
V\rangle^{2}_{0}\right)=K_{21}s_{a}s_{b}+K_{22}\left(s^{2}_{a}+s^{2}_{b}\right)
+K_{23}\left(s^{4}_{a}+s^{4}_{b}\right),\label{H'2}
\end{eqnarray}
with
\begin{subequations}\begin{eqnarray}
K_{21}&=&\frac{105K^{2}u^{2}}{4b^{5}},\\
K_{22}&=&\frac{105K^{2}u^{2}}{8b^{5}},\\
K_{23}&=&\frac{87K^{4}u^{2}}{16b^{6}}.
\end{eqnarray}\label{20}\end{subequations}
Combining Eqs.(\ref{00}), (\ref{10}) and (\ref{20}), we obtain the
effective Hamiltonian after RG transformation
\begin{eqnarray}
H'_{\textrm{eff}}&=&(K_{01}+K_{11}+K_{21})s_{a}s_{b}+(K_{02}+K_{12}+K_{22})\left(s^{2}_{a}+s^{2}_{b}\right)\\\nonumber
&+&(K_{03}+K_{13}+K_{23})\left(s^{4}_{a}+s^{4}_{b}\right),\label{H'eff}
\end{eqnarray}
In order to get the recursion relations of the RG transformation, we
rescale the spins by $s'_{i}=\xi_{i}s_{i}$ ($i=a,b$) and then Eq.
(\ref{H'eff}) can be rewritten as
\begin{eqnarray}
H'_{\textrm{eff}}=K's'_{a}s'_{b}-\frac{b}{2}\left(s'^{2}_{a}+s'^{2}_{b}\right)
-u'\left(s'^{4}_{a}+s'^{4}_{b}\right),\label{effH'}
\end{eqnarray}
where
\begin{subequations}\begin{eqnarray}
\xi_{i}&=&\sqrt{\frac{-2(K_{11}+K_{12}+K_{13})}{b}},\\
K'&=&\frac{(K_{01}+K_{02}+K_{03})}{\xi^{2}_{i}},\\
u'&=&-\frac{(K_{21}+K_{22}+K_{23})}{\xi^{4}_{i}}.
\end{eqnarray}\label{rscale}\end{subequations}
Setting $K'=K=K^{*}$ and $u'=u=u^{*}$, the fixed points of the
system can be found are
\begin{subequations}\begin{eqnarray}
A:\hspace{0.3cm} K^{*}&=&0, u^{*}=0;\\
B:\hspace{0.3cm} K^{*}&=&b, u^{*}=0;\\
C:\hspace{0.3cm} K^{*}&=&0.436b, u^{*}=0.370b^{2};
\end{eqnarray}\end{subequations}\label{20}
in which $A$ is found to be a nonphysical steady fixed point while
$B$ is a Gaussian fixed point and $C$ a Wilson-Fisher one.

Actually, in the three fixed point, the Wilson-Fisher one is most
crucial for the critical properties of the $S^{4}$ system. By
expanding $K'$ and $u'$ in the neighborhood of the Wilson-Fisher
fixed point $C$ and preserving only the linear items, the
transformation matrix is obtained to be
\begin{eqnarray}
R_{\textrm{L}}(K,u)=\left(
\begin{array}{cc}
\frac{\partial K'}{\partial K} & \frac{\partial K'}{\partial
u} \\
\frac{\partial u'}{\partial K} & \frac{\partial u'}{\partial
u} \\
\end{array}
\right)_{C}=\left(
\begin{array}{cc}
2.855 & 2.895 \\
1.356 & 2.447 \\
\end{array}
\right)
\label{matrixa},
\end{eqnarray}
with two eigenvalues $\lambda_{1}=4.643>1$ and $\lambda_{2}=0.659<1$
revealing that $C$ is a barrier point. Thus a scale power can be
obtained according to the scaling theory
$p=\frac{\textrm{ln}\lambda_{1}}{d\textrm{ln}L}=1.107$ where $d=2$
is the dimensionality of the triangular lattice and $L=2$ is the
scaling factor. The correlation length exponent can be obtained from
the relations $\nu=\frac{1}{pd}$ to be $\nu=0.451$ smaller than the
Gaussian case which contains only two-body interactions
\cite{chyw3,sliyang} while in good conformity with previous $S^{4}$
system studies on other lattices \cite{liy,yin}.

%

\subsection{external field}

In order to obtain all the critical exponents, we study in the
following the case of $S^{4}$ system in an external field $h$ where
the effective Hamiltonian can be written as Eq.(\ref{eH}). For
simplicity, we give out here only the primary results while leaving
the process of derivation in the appendixes. Considering only
nearest-neighbor interactions, three fixed points can also be
obtained by following similar RG transformation procedures. They are
\begin{subequations}\begin{eqnarray}
A:\hspace{0.3cm} K^{*}&=&0, u^{*}=0, h^{*}=0;\\
B:\hspace{0.3cm} K^{*}&=&b, u^{*}=0, h^{*}=0;\\
C:\hspace{0.3cm} K^{*}&=&0.436b, u^{*}=0.352b^{2}, h^{*}=0;
\end{eqnarray}\label{20}\end{subequations}
the transformation matrix in the neighborhood of the Wilson-Fisher
fixed point $C$ is now
\begin{eqnarray}
R_{\textrm{L}}(K,u)=\left(
\begin{array}{ccc}
\frac{\partial K'}{\partial K} & \frac{\partial K'}{\partial u} &
\frac{\partial K'}{\partial
h} \\
\frac{\partial u'}{\partial K} & \frac{\partial u'}{\partial u} &
\frac{\partial u'}{\partial
h} \\
\frac{\partial h'}{\partial K} & \frac{\partial h'}{\partial u} &
\frac{\partial h'}{\partial
h} \\
\end{array}
\right)_{C}=\left(
\begin{array}{ccc}
2.872 & 3.032 & 0 \\
1.307 & 2.450 & 0 \\
0     & 0     & 2.875 \\
\end{array}
\right) \label{matrixa},
\end{eqnarray}
with three eigenvalues $\lambda_{1}=4.663>1$, $\lambda_{2}=0.659<1$
and $\lambda_{3}=2.875>1$. Where $\lambda_{1}$ and $\lambda_{2}$ is
related to the temperature while $\lambda_{3}$ is related to the
external field. Then we can obtain the two scale powers according to
the scaling theory
\begin{eqnarray}
p=\frac{\textrm{ln}\lambda_{1}}{d\textrm{ln}L}=1.110, \hspace{1cm}
q=\frac{\textrm{ln}\lambda_{3}}{d\textrm{ln}L}=0.762;
\end{eqnarray}
The critical exponents at this point can be obtained following
$\alpha=\frac{2p-1}{p}$, $\beta=\frac{1-q}{p}$,
$\gamma=\frac{2q-1}{p}$, $\delta=\frac{q}{1-q}$, $\eta=2+d(1-2q)$,
$\nu=\frac{1}{pd}$ to be $\alpha=1.099$, $\beta=0.215$,
$\gamma=0.471$, $\delta=3.197$, $\eta=0.953$, $\nu=0.450$ describing
in the neighborhood of the critical point the variations in heat,
magnetization, magnetic field, magnetic susceptibility, correlation
function and correlation length, respectively.

\subsection{next-nearest neighbor interaction} \label{sbsec32}

\begin{figure}
\centering
\includegraphics[scale=1.2]{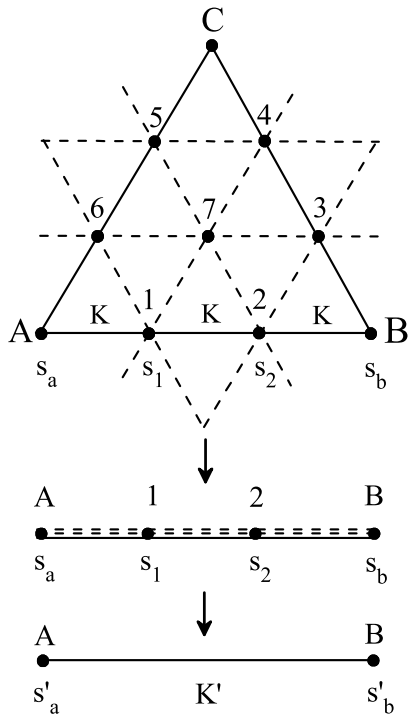}
\caption{Bond-moving and decimation procedures for the renormalized
bond $K'$ in the case of considering next nearest-neighbor
interactions.\label{pro3}}
\end{figure}

For the enriching of the $S^{4}$ results, we present in this paper
also the case of considering next-nearest neighbor interactions. All
the approximation of the effective Hamiltonian after RG
transformation can be obtained also by combining the generalized
bond-moving operations and the cumulative expansion techniques. The
only change need to be made is to select a cluster of containing
several sites (as is shown in Fig.\ref{pro3}) as a basic unit for
recursion. The various coefficients of $H'_{\textrm{eff}}$ after RG
derivation is also left in the appendixes for simplicity, we give
out here only the primary results.

In the case of considering next-nearest neighbor interactions, we
can also obtain by RG transformation three fixed points
\begin{subequations}\begin{eqnarray}
A:\hspace{0.3cm} K^{*}&=&0, u^{*}=0;\\
B:\hspace{0.3cm} K^{*}&=&0.829b, u^{*}=0;\\
C:\hspace{0.3cm} K^{*}&=&0.257b, u^{*}=5.405b^{2};
\end{eqnarray}\label{20}\end{subequations}
in which $A$ is found to be a nonphysical steady fixed point while
$B$ is a Gaussian fixed point and $C$ a Wilson-Fisher one. After
some algebra we can obtain the correlation length critical exponent
at the Wilson-Fisher fixed point to be $\nu=0.341$ also in
conformity with previous studies. However we find it becomes
relatively smaller than previous results. This reveals the effects
of next nearest neighbor interactions.

\section{summary and discussions} \label{sec4}
In summary, we have presented in this paper a detailed study on the
critical behavior of the $S^{4}$ system constructed on the
triangular lattices by combining the generalized bond-moving RG
transformation method and the cumulative expansion techniques. In
all cases it is found to have three fixed points. The critical
exponents obtained near the Wilson-Fisher fixed points are found
relatively smaller than that of the Gaussian model and different
from those of the Ising model. This reveals the particular effect of
the four-body interactions on the critical properties of the $S^{4}$
system. In further by comparing all the results of correlation
length exponent obtained in all cases we find that it becomes
smaller and smaller with the increasing of the complexity of the
system. This reveals the decisive influence of the system complexity
on the critical behavior.

\section * {ACKNOWLEDGEMENTS}

This work was supported by the Scientific Research Starting
Foundation (Grant No.Bsqd2008053), Youth Foundation (Grant
No.XJ201009) and the Scientific Research Training Plan for
Undergraduate Students (Grant No.2010A023) of Qufu Normal
University.

\appendix
\section {RG coefficients of $H'_{\textrm{eff}}$ in external field}
For simplicity we give out here only the various RG coefficients of
$H'_{\textrm{eff}}$ from zero to second order approximations instead
of showing all the details of derivation that is very similar to
those have appeared in this paper hereinbefore. In the case of
considering the influence of an external field $h$, these
coefficients are
\begin{eqnarray}
K_{01}&=&\frac{K^{2}}{b},\\
K_{02}&=&\frac{K^{2}-2b^{2}}{2b},\\
K_{03}&=&-2u,\\
K_{04}&=&\frac{(K-2b)h}{b},\\
K_{11}&=&-\frac{3K^{2}u}{b^{4}}(b+h^{2}),\\
K_{12}&=&-\frac{3K^{2}u}{2b^{4}}(b+h^{2}),\\
K_{13}&=&-\frac{K^{4}u}{4b^{4}},\\
K_{14}&=&-\frac{Kuh}{b^{4}}(3b+h^{2}),\\
K_{21}&=&\frac{3K^{2}u^{2}}{4b^{5}}(35+\frac{93h^{2}}{b}+\frac{9h^{4}}{b^{2}}+\frac{h^{6}}{b^{3}}),\\
K_{22}&=&\frac{3K^{2}u^{2}}{8b^{5}}(35+\frac{93h^{2}}{b}+\frac{9h^{4}}{b^{2}}+\frac{h^{6}}{b^{3}}),\\
K_{23}&=&\frac{K^{4}u^{2}}{16b^{6}}(87+\frac{174h^{2}}{b}+\frac{17h^{4}}{b^{2}}),\\
K_{24}&=&\frac{hKu^{2}}{4b^{5}}(105+\frac{105h^{2}}{b}+\frac{21h^{4}}{b^{2}}+\frac{h^{6}}{b^{3}}).
\end{eqnarray}\label{a}where $K_{i4}$ $(i=0,1,2)$ represents the various RG coefficients of
$H'_{\textrm{eff}}$ which is related with the external field.
\section {RG coefficients of next-nearest $H'_{\textrm{eff}}$}
The RG coefficients for various order approximation of the effective
Hamiltonian in the case of considering next-nearest neighbor
interactions are
\begin{eqnarray}
K_{01}&=&\frac{48K^{2}b+51K^{3}}{16(4b^{2}-K^{2})},\\
K_{02}&=&-\frac{96b^{3}-75K^{2}b-12K^{2}}{16(4b^{2}-K^{2})},\\
K_{03}&=&-3u,\\
K_{11}&=&\frac{3K^{2}u}{4b^{3}(4b^{2}-K^{2})^{2}},\\
K_{12}&=&-\frac{K^{2}u(272b^{4}+128kb^{3}+136K^{2}b^{2}-17K^{4})}{64b^{3}(4b^{2}-K^{2})^{2}},\\
K_{13}&=&-\frac{K^{4}u(4112b^{4}+2176Kb^{3}+2824K^{2}b^{2}+1088K^{3}b+771K^{4})}{4096b^{4}(4b^{2}-K^{2})^{2}},\\
K_{21}&=&-\frac{2K^{2}u^{2}(656b^{5}+1802b^{4}K+452b^{3}K^{2}+374b^{2}K^{3}+56bK^{4}+17K^{5})}{9(4b^{2}-K^{2})^{5}},\\
K_{22}&=&-\frac{K^{2}u^{2}(2584b^{5}+1408b^{4}K+1258b^{3}K^{2}+176b^{2}K^{3}+119bK^{4}+8K^{5})}{9(4b^{2}-K^{2})^{5}},\\\nonumber
K_{23}&=&\frac{K^{4}u^{2}}{24576(4b^{3}-bK^{2})^{6}}(15790080b^{12}+9469952b^{11}K+6037504b^{10}K^{2}\\
\nonumber&&+5431296b^{9}K^{3}+4420864b^{8}K^{4}+1114112b^{7}K^{5}+185088b^{6}K^{6}\\&&-208896b^{5}K^{7}
-152592b^{4}K^{8}+13056b^{3}K^{9}+20808b^{2}K^{10}-867K^{12}).
\end{eqnarray}\label{b}where the next nearest-neighbor interaction is set one quarter in
strength as the nearest-neighbor interaction.


\begin{thebibliography}{References}

\bibitem{msrg} L. E. Reichl, A Modern Course in Statistical Physics, second ed., New York, (1998).

\bibitem{liy} Y. Li, X. M. Kong, Physica A \textbf{356}, 589 (2005).

\bibitem{yin} X. C. Yin, H. Yin and X. M. Kong, Acta. Phys. Sin. \textbf{55}, 4901 (2006).

\bibitem{chyw3} C. Y. Wang, W. X. Yang, Z. W. Yan, et al. Commun. Theor. Phys. \textbf{57}, 717 (2012).

\bibitem{liu} J. Liu, X. M. Kong, Y. P. Li and J. Y. Huang, Acta. Phys. Sin. \textbf{53}, 2275 (2004).

\bibitem{s4} S. Mac, Modern Theory of Critical Phenomena, New York, 1976.

\bibitem{Genfen} Y. Gefen, A. Aharony and B. B. Mandelbrot, J. Phys. A: Math. Gen. \textbf{16} 1267
(1983).
\bibitem{Genfen2} Y. Gefen, A. Aharony and B. B. Mandelbrot, J. Phys. A: Math. Gen. \textbf{17}
435 (1984).
\bibitem{Genfen3} Y. Gefen, A. Aharony and B. B. Mandelbrot, J. Phys. A: Math. Gen. \textbf{17} 1277
(1984).
\bibitem{sliyang} S. Li and Z. R. Yang, Phys. Rev, E \textbf{55} 6656 (l997).

\end{thebibliography}
\end{document}